\documentclass[12pt]{article}
\usepackage{amsmath}
\usepackage{amssymb}
\usepackage{amsthm}

\theoremstyle{plain}
\newtheorem{theorem}{Theorem}
\newtheorem{proposition}{Proposition}

\theoremstyle{definition}
\newtheorem{definition}{Definition}
\newtheorem{algorithm}{ALGORITHM}

\theoremstyle{remark}
\newtheorem{remark}{Remark}

\newcommand{\Ind}{\mathop{\rm Ind}\nolimits}
\newcommand{\bst}{\mathop{\rm bst}\nolimits}
\newcommand{\rank}{\mathop{\rm rank}\nolimits}

\begin{document}

\title{Index Function and Minimal Cycles}
\author{A.V. Lapteva, E.I. Yakovlev\\
{\it Nizhni Novgorod State University}}
\date{}
\maketitle

\thispagestyle{empty}

\begin{abstract}
Let $P$ be a closed triangulated manifold, $\dim{P}=n$. We
consider the group of simplicial 1-chains $C_1(P)=
C_1(P,\Bbb{Z}_2)$ and the homology group $H_1(P)=
H_1(P,\Bbb{Z}_2)$. We also use some nonnegative weighting function
$L:C_1(P)\to\Bbb{R}$. For any homological class $[x]\in{H}_1(P)$
method proposed in article builds a cycle $z\in[x]$ with minimal
weight $L(z)$. The main idea is in using a simplicial scheme of
space of the regular covering $p:\hat{P}\to{P}$ with automorphism
group $G\cong{H}_1(P)$. We construct this covering applying index
function $J:C_1(P)\to\Bbb{Z}_2^r$ relative to any basis of group
$H_{n-1}(P)$, $r= \rank{H_{n-1}(P)}$.

\end{abstract}

\section{Index Function}

Consider a triangulated closed manifold $P$, $\dim{P}=n$, and a
basis $[z_1^{n-1}],\dots,[z_r^{n-1}]$ of homology group
$H_{n-1}(P)=H_{n-1}(P,\Bbb{Z}_2)$. Let
$\Ind:H_1(P)\times{H}_{n-1}(P)\to\Bbb{Z}_2$ be intersection index.

\begin{definition}\label{Index function}
Define the homomorphism $J_0: Z_1(P)\to\Bbb{Z}_2^r$ by the
formulas $J_0^k(y)=\Ind([y],[z_k^{n-1}])$, $k=1,\dots,r$,
$J_0=(J_0^1,\dots,J_0^r)$. We call its arbitrary extension
$J:C_1(P)\to\Bbb{Z}_2^r$ index function. For any chain
$x\in{C}_1(P)$ the value $J(x)$ is called its index relative to
the basis $[z_1^{n-1}],\dots,[z_r^{n-1}]$.
\end{definition}

\begin{remark}
Index function $J:C_1(P)\to\Bbb{Z}_2^r$ is not uniquely defined,
however we can use this function for solve our problems.
\end{remark}

\begin{proposition}\label{Index function properties}
If $J:C_1(P)\to\Bbb{Z}_2^r$ is index function relative to the
basis $\left\{[z_1^{n-1}],\dots,[z_r^{n-1}]\right\}$ of group
$H_{n-1}(P)$, $x,y\in{C}_1(P)$ and $\partial{x}=\partial{y}$, then
$J(x)=J(y)$ if and only if $x\sim{y}$.
\end{proposition}

\begin{proof}
Let $\left\{[z_1^1],\dots,[z_r^1]\right\}$ be a basis of group
$H_1(P)=H_1(P,\Bbb{Z}_2)$, that is dual to the given basis
$\left\{[z_1^{n-1}],\dots,[z_r^{n-1}]\right\}$. Assume now that
$z=x+y$. Then $z\in{Z}_1(P)$ and $[z]=\sum_{i=1}^r{l}^i[z_i^1]$,
where $l_i\in\Bbb{Z}_2$. This implies that
$J^k(z)=\Ind([z],[z_k^{n-1}])=l^k$ for all $k=1,\dots,r$. So
$J(x)=J(y)$ if and only if $l^1=\dots=l^r=0$. And this latter
expression is equivalent to equality $[z]=0$.
\end{proof}

\begin{algorithm}\label{Construction index function}
{\bf Construction of index function relative to the basis of group
$H_{n-1}(P)$.}
\smallskip

\noindent {\bf Input:}

1) simple basis cycles $z_1^{n-1},z_2^{n-1},\dots,z_r^{n-1}$ which
are lists of $(n-1)$-dimensional simplices;

2) list $K^1(P)$ of edges for polyhedron $P$;

3) lists $K_1^n(P,z_1^{n-1}),\dots,K_r^n(P,z_r^{n-1})$ consisting
of $n$-dimensional simplices from neighbourhoods of cycles
$z_1^{n-1},\dots,z_r^{n-1}$ respectively;
\smallskip

\noindent {\bf Output:}

1) vectors $J(a)=\left(J^1(a),\dots,J^r(a)\right)\in\Bbb{Z}_2^r$
for all edges $a\in K^1(P)$;

2) chains $M_1,\dots,M_r$ of edges indexed relative to cycles
$z_1^{n-1},\dots,z_r^{n-1}$ respectively;

3) lists $M_k(u)$, $k=1,\dots,r$ of edges,
that we add to $M_k$ when considering vertex $u$ of cycle $z_k^{n-1}$;

4) sets $\Sigma_k(u)$, $k=1,\dots,r$ of $n$-simplices
incident to edges from $M_k(u)$.
\smallskip

\noindent{\bf Algorithm Description.}
\smallskip

{\bf Step 0.} For all $k=1,\dots,r$ execute steps 1 -- 3.

{\bf Step 1. Start operations.} Assume $M_k=\emptyset$,
$J^k(a):=0$ for all $a\in K^1(P)$. We denote $z_k^{n-1}$ by $X$
and $K_k^n(P,z_k^ {n-1})$ by $K^n(P,X)$. We create then lists of
vertices and edges for all simplices of cycle $X$, $K^0(X)$ and
$K^1(X)$ respectively.

{\bf Step 2. Indexing edges that do not belong to the cycle.} For
each vertex $u\in K^0(X)$ execute steps 2.1 -- 2.4.

{\bf Step 2.1. Initializing vertex neighbourhood.} Create a list
$K^n(P,u)\subset K^n(P,X)$ of $n$-dimensional simplices of the
polyhedron $P$, that contain $u$, and a list $K^{n-1}(P,u)$ of all
$(n-1)$-dimensional faces of simplices from $K^n(P,u)$. At the
same time, for each simplex $\sigma^{n-1}\in K^{n-1}(P,u)$ we get
a list $\partial^{n-1}(\sigma^{n-1},u)$ of $n$-dimensional
simplices from $K^n(P,u)$ those are incident to $\sigma^{n-1}$,
and assume $\mu(\sigma^{n-1}):=0$. Then we create empty lists
$M_k(u):=\emptyset$ and $\Sigma_k(u):=\emptyset$.

{\bf Step 2.2. Creating the queue to keep $n$-simplices.} We chose
a simplex $\sigma_0^n\in K^n(P,u)$, create a queue
$R:=\{\sigma_0^n\}$ and remove $\sigma_0^n$ from $K^n(P,u)$.

{\bf Step 2.3. Main procedure of the Algorithm.} While the queue
$R$ is not empty we will do the following actions. Take the first
simplex $\sigma^n\in R$ and remove it from the queue $R$. For each
$(n-1)$-dimensional face $\sigma^{n-1}$ of the simplex $\sigma^n$
we check the following: whether it belongs to the cycle $X$,
whether $\mu(\sigma^{n-1})$ is equal to zero, whether the list
$\partial^{n-1}(\sigma^{n-1},u)$ contains any simplices different
from $\sigma^n$. If all above conditions are satisfied we will
execute steps 2.3.1 -- 2.3.2.

{\bf Step 2.3.1.} Take the simplex
$\sigma_*^n\in\partial^{n-1}(\sigma^
{n-1},u)\setminus\{\sigma^n\}$, remove it from $K^n(P,u)$ and
enqueue to $R$; set $\mu(\sigma^{n-1}):=1$ and $\Sigma_k(u)=
\Sigma_k(u)\cup\{\sigma_*^n\}$.

{\bf Step 2.3.2.} For all vertices $w\ne u$ of the simplex
$\sigma^{n-1}$ we check whether the edge $a=[uw]$ is in the list
$K^1(X)$; having $a\notin K^1(X)$ set $J^k(a):=J^k(a)+1\mod 2$,
$M_k(u)=M_k(u)\cup\{a\}$, $M_k:=M_k+a\mod 2$.

{\bf Step 2.4. Main procedure repeated.} If the list $M_k(u)$ is
empty then go back to step 2.2.

{\bf Step 3. Indexing the edges of cycle.} For each edge
$a=[uv]\in K^1(X)$ we search any edges $b\in M_k(u)$ and $c\in
M_k(v)$ such that $b\cap c\ne\emptyset$ and that $a$, $b$ and $c$
are sides of some triangle of polyhedron $P$. If the edges $b$ and
$c$ do not exist then we set $J^k(a):=1$ and $M_k=M_k+a\mod 2$.

\smallskip

\noindent{\bf End of algorithm.}
\end{algorithm}

\begin{theorem}\label{th of index function construction}
If $P$ is a closed $n$-dimensional manifold,
$z_1^{n-1},\dots,z_r^{n-1}$ are simple cycles, $x=a_1+\dots+a_l\in
C_1(P)$ and $J(x)=\sum_{i=1}^l {J(a_i)}$, then the vector
$J(x)=\left( J^1(x),\dots,J^r(x)\right)\in\Bbb{Z}_2^r$ is index of
the chain $x\in C_1(P)$ relative to the basis
$[z_1^{n-1}],\dots,[z_r^{n-1}]$ of group $H_{n-1}(P)$.
\end{theorem}

\begin{proof}
Let $x\in Z_1(P)$. We will prove that
$J^k(x)=\Ind([x],[z_k^{n-1}])$ for all $k=1,\dots,r$.

Set $z_0^*=z_k^{n-1}$. For all $p=1,\dots,N$ we will make the
following constructions; here $N$ is power of the set
$K^0(z_k^{n-1})$.

Consider vertex $u_p\in K^0(z_k^{n-1})$ and its barycentric star
$\bst(u_p,P)$.

Let $\Sigma_k^*(u_p)$ be the set of all $n$-simpleces from the
barycentric subdivision of $\Sigma_k(u_p)$. Construct the chain
$c(u_p)$ of simpleces $\sigma_1\in \bst(u_p,P)\cap
\Sigma_k^*(u_p)$.

Then we write the chain boundary $c(u_p)$ as a sum $Y_1+Y_2$,
where $Y_1$ is the sum of all its $(n-1)$-dimensional simplices,
that belong to cycle $z_k^{n-1}$ and $Y_2$ is the sum of all
remaining simplices from the chain $\partial c(u_p)$. Set
$z_p^*=z_{p-1}^*+Y_1+Y_2\mod 2$.

By construction $z_p^* \sim z_{p-1}^*$ for all $p=1,\dots,N$.
Hence, the cycle $z^*=z_N^*$ is homologous to the cycle
$z_k^{n-1}=z_0^*$.

Let now prove that for any edge $a=[uv]\in K^1(P)$ and
$\sigma_b\in\bst(a)$ the simplex $\sigma_b$ belong to $z^*$ if and
only if $a\in M_k$.

Let view all possible positions of the edge $a$. At the same time
we also agree to think that ${M}_k(u)=\emptyset$ and that
$\Sigma_k(u)=\emptyset$ for all $u\notin{K}^0(z_k^{n-1})$.

0. If $a\notin{M}_k(u)\cup{M}_k(v)$ and $a\notin{K}^1(z_k^{n-1})$,
then according to the algorithm $a\notin{M}_k$. On the other hand,
the edge $a$ can not be incident to simplices from the lists
$\Sigma_k(u)$ and $\Sigma_k(v)$ and hence $\sigma_b\notin z^*$.

1. Let $u\in K^0(z_k^{n-1})$, $a\in M_k(u)$ and $v\notin
K^0(z_k^{n-1})$. Then the edge $a$ will be still in the chain
$M_k$ when algorithm \ref{Construction index function} is
completed. At the same time the barycentric star $\bst(a)$ belongs
to the boundary of the chain $c(u)$ and does not belong to the
cycle $z_k^{n-1}$. Thus in this case $a\in{M}_k$ and the chain
$\bst(a)$ belongs to the cycle $z^*$.

2. Further, assume that $u,v\in K^0(z_k^{n-1})$ and $a\in M_k(u)$.
At that, $a\notin K^1(z_k^{n-1})$.

2.1. If $a\in M_k(v)$, then $a\notin M_k$, and simplices of its
barycentric star will be added twice to the initial cycle
$z_k^{n-1}$ and will not be in the resulting cycle $z^*$.

2.2. If $a\notin M_k(v)$, then $a\in M_k$ and any simplex
$\sigma_b\in\bst(a)$ is added to the cycle $z^*$ exactly once. So
$\sigma_b\in z^*$.

3. Finally, let $a\in K^1(z_k^{n-1})$.

3.1. Let assume that the condition from step 3 of algorithm
\ref{Construction index function} is satisfied, i.e.:
\begin{itemize}
\item[($*$)] there exist edges $b\in M_k(u)$ and $c\in M_k(v)$
such that $b\cap c\ne\emptyset$ and that $a$, $b$ and $c$ are
sides of some triangle $\sigma'\in{K}^2(P)$.
\end{itemize}

In this case, according to the algorithm $a\notin M_k$.

Let view all triangles $\sigma'$ from ($*$), and all
$n$-dimensional simplices incident to them. The such $n$-simplices
belong both to $\Sigma_k(u)$ and $\Sigma_k(v)$. Consider
$n$-dimensional simplex $\sigma$, $\sigma_b\in\sigma$. If
$\sigma_b\in\bst(a)$, then $\sigma$ either belong to the both sets
$\Sigma_k(u)$ and $\Sigma_k(v)$ or does not belong to them. Hence,
the simplex $\sigma_b$ either is not added to the cycle $z^*$ or
is added twice. Therefore $\sigma_b\notin z^*$.

3.2. Assume now that condition ($*$) is not satisfied. Then
according to step 3 of algorithm \ref{Construction index
function}, $a\in M_k$.

Barycentric star $\bst(a)$ of the edge $a=[uv]$ belongs to
the union $D(a)$ of all $n$-simplices that contain the edge $a$.
We will prove that the sub-polyhedron $D(a)$ belongs to the union of simplices from
the sets $\Sigma_k(u)$ and $\Sigma_k(v)$.

Cycle $z_k^{n-1}$ divides $D(a)$ into two components of strong
connectivity $D^+(a)$ and $D^-(a)$.

By construction the set $\Sigma_k(u)$ can not be empty. Moreover,
if the simplex $\sigma~\in z_k^{n-1}$ is incident to the vertex
$u$, then $\sigma~$ is a face of some $n$-simplex from
$\Sigma_k(u)$. So there exists a simplex $\sigma^n\in\Sigma_k(u)$
that contains the edge $a$.

Let the simplex $\sigma^n$ belongs to $D^+(a)$. Then under the
strong connectivity $D^+(a)$ and according to algorithm
\ref{Construction index function}, all $n$-simplices from $D^+(a)$
also belong to $\Sigma_k(u)$.

This implies, in accordance with our assumption, that no
$n$-simplex from $D^+(a)$ can belong to the set $\Sigma_k(v)$.

The set $\Sigma_k(v)$ can not be empty also. Since each simplex of
$z_k^{n-1}$ incident to the vertex $v$ is a face of some
$n$-simplex from $\Sigma_k(v)$, it follows that there exists a
simplex $\sigma^n_*\in\Sigma_k(v)$ that contains the edge $a$. By
the above proof $\sigma^n_*$ belongs to $D^-(a)$. Then all
$n$-simplices from $D^-(a)$ belong to the set $\Sigma_k(v)$ too.
Consequently all $n$-simplices of the polyhedron $D(a)=D^+(a)\cup
D^-(a)$ belong either to the set $\Sigma_k(u)$ or to
$\Sigma_k(v)$.

Consider $\sigma_b\in\bst(a)$. If there exists a simplex
$\sigma\in\Sigma_k(u)$ containing $\sigma_b$, then $\sigma\notin
\Sigma_k(v)$. Otherwise, in accordance to the above proof, there
is a simplex $\tilde{\sigma}\in\Sigma_k(v)$ such that
$\sigma_b\subset\tilde{\sigma}$. It follows that $\sigma_b$ is
involved in the cycle $z^*$ exactly once, so $\sigma_b\in z^*$.

Thus we have proved that the cycle $z^*\sim z^{n-1}_k$ consists of
barycentric stars of the edges from chain $M_k$. That means that
this cycle intersects transversally only the edges of the cycle
$x$, that are in the list $M_k$. According to algorithm
\ref{Construction index function} $J^k(a)=1$ for all $a\in M_k$
and $J^k(b)=0$ for all edges $b\notin M_k$. So
$$\Ind([x],[z_k^{n-1}])=\Ind([x],[z^*])=\sum\limits_{a\in x} J^k(a)\mod2=J^k(x).$$
\end{proof}

\begin{remark}
The fact that $[z_1^{n-1}],\dots,[z_r^{n-1}]$ is a basis of group
$H_{n-1}(P)$ has no impact on the behaviour of algorithm
\ref{Construction index function}. So we can apply this algorithm
to an arbitrary set of simple $(n-1)$-dimensional cycles of the
manifold $P$. In particular this set may consist of only one cycle
$z^{n-1}$. Then we will get a function $J:C_1(P)\to\Bbb{Z}_2$ such
that $\sum_{i=1}^l {J(a_i)}=\Ind([x],[z^{n-1}])$ for
$x=a_1+\dots+a_l\in{Z}_1(P)$. So we can use algorithm
\ref{Construction index function} to compute the intersection
index of a given $(n-1)$-cycle $z^{n-1}\in{Z}_{n-1}(P)$ with any
one-dimensional cycle of the manifold $P$.
\end{remark}

\begin{remark}
We can find any basis $[z_1^{n-1}],\dots,[z_r^{n-1}]$ of group
$H_{n-1}(P)$ using standard matrix algorithm (see, for example,
\cite{ZT}). If $n=2$, we also can apply algorithms that don't use
incidence matrices (see \cite{V_Y,YGR}).
\end{remark}

\section{Regular Covering with the Automorphism Group $H_1(P)$}

Let $P$ be a $n$-dimensional triangulated closed manifold and
$S=(V,K)$ be its simplicial scheme. We will construct an abstract
simplicial scheme $\hat{S}=(\hat{V},\hat{K})$ as follows.

Set $\hat{V}=V\times{G}$, where $G=\Bbb{Z}_2^r$. Let
$\hat{v}_0,\hat{v}_1,\dots,\hat{v}_m\in\hat{V}$, where
$\hat{v}_i=(v_i,b_i)$ for all $i=0,1,\dots,m$. We will think that
$\{\hat{v}_0,\hat{v}_1,\dots,\hat{v}_m\}\in\hat{K}$ if the below
conditions are satisfied:
\begin{itemize}
\item[(U1)] $\{v_0,v_1,\dots,v_m\}\in{K}$;
\item[(U2)] $g_0+g_i=J([v_0v_i])$ for any $i=1,\dots,m$;
here $J([v_0v_i])$ is the index of the edge $[v_0v_i]$.
\end{itemize}

\begin{remark}\label{0ij}
When the conditions (U1) and (U2) are satisfied the equalities
$g_i+g_j=J([v_iv_j])$ are also true for all $i,j=1,\dots,m$. In
fact, according to (U1), the cycle $z=[v_jv_i]+[v_iv_0]+[v_0v_j]$
is homologous to zero. So $J([v_iv_j])=J([v_iv_0])+J([v_0v_j])$.
By invoking (U2) we can have these equalities
$J([v_iv_j])=g_i+g_0+g_0+g_j=g_i+g_j$.
\end{remark}

Let define now a mapping $p^{\,0}:\hat{V}\to{V}$ and a left action
$\lambda^0:G\times\hat{V}\to\hat{V}$ of group $G$ on $\hat{V}$,
assuming
\begin{equation}\label{def p and G action}
p^{\,0}((v,g))=v\quad {\text{and}}\quad \lambda^0(g',(v,g))=g'\cdot(v,g)=(v,g'+g)
\end{equation}
for all $(v,g)\in\hat{V}$ and $g'\in{G}$.

Let $\hat{P}$ define some realization of the scheme
$\hat{S}=(\hat{V},\hat{K})$. At that we identify the set of
vertices of the polyhedron $\hat{P}$ with $\hat{V}$.

\begin{proposition}\label{simplicial covering with group G}
For the mapping $p^{\,0}:\hat{V}\to{V}$ there
exists the unique continuation $p:\hat{P}\to{P}$ that is
simplicial regular covering with a group of covering
transformations $G\cong{H}_1(P)$.
\end{proposition}

\begin{proof}
Simplicial and surjective properties of the mapping $p^{\,0}$
follow directly from its definition and from the construction of
the complex $\hat{K}$. If
$\hat{s}=\{(v_0,g_0),(v_1,g_1),\dots,(v_m,g_m)\}\in\hat{K}$, then
$\{v_0,v_1,\dots,v_m\}\in{K}$ and $g_0+g_i=J([v_0v_i])$ for all
$i=1,\dots,m$. On the other hand,
$g\cdot\hat{s}=\{(v_0,g+g_0),(v_1,g+g_1),\dots,(v_m,g+g_m)\}$ for
an arbitrary $g\in{G}$. Since $g+g_0+g+g_i=g_0+g_i=J([v_0v_i])$,
then $g\cdot\hat{s}\in\hat{K}$. So the action $\lambda^0$ is also
simplicial.

Let $s=\{v_0,v_1,\dots,v_m\}\in{K}$ and
$\hat{v}_0\in({p}^0)^{-1}(v_0)$. Then $\hat{v}_0=(v_0,g_0)$, where
$g_0\in{G}$. Set $g_i=g_0+J([v_0v_i])$ and $\hat{v}_i=(v_i,g_i)$
for all $i=1,\dots,m$. At that
$\hat{s}=\{\hat{v}_0,\hat{v}_1,\dots,\hat{v}_m\}\in\hat{K}$,
$\hat{v}_0\in\hat{s}$ and $p^{\,0}(\hat{s})=s$. Hence, the mapping
$p^{\,0}$ has the following property:
\begin{itemize}
\item[\rm(C1)] for each abstract simplex $s\in{K}$
and for any vertex $\hat{v}\in({p}^0)^{-1}(s)$ there is the unique
abstract simplex $\hat{s}\in\hat{K}$ containing the vertex
$\hat{v}$ and satisfying the equality $p^{\,0}(\hat{s})=s$.
\end{itemize}

Let choose an abstract simplex
$\hat{s}=\{\hat{v}_0,\hat{v}_1,\dots,\hat{v}_m\}\in\hat{K}$, and
an element $g$ of group $G$ and assume that
$g\cdot\hat{s}=\hat{s}$. Then $\hat{v}_i=(v_i,g_i)$ and
$g\cdot\hat{v}_i=(v_i,g+g_i)$ for all $i=1,\dots,m$. At the same
time it follows from the equality $g\cdot\hat{s}=\hat{s}$ that
$(v_0,g+g_0)=(v_k,g_k)$ for some $k\in\{0,1,\dots,m\}$. The latter
is possible only if $k=0$ and $g=0$. Thus the action $\lambda^0$
has the following property:
\begin{itemize}
\item[\rm(C2)] if $g\cdot\hat{s}=\hat{s}$
for at least one non-empty simplex $\hat{s}\in\hat{K}$, then $g$
is the neutral element of the group $G$.
\end{itemize}

Let now consider the simplices
$\hat{s}=\{(v_0,g_0),(v_1,g_1),\dots,(v_m,g_m)\}$ and $\hat{s}'$
of the complex $\hat{K}$.

First, if $g\in{G}$ and $\hat{s}'=g\cdot\hat{s}$, then
$\hat{s}'=\{(v_0,g+g_0),(v_1,g+g_1),\dots,(v_m,g+g_m)\}$. At that
$p^{\,0}(\hat{s}')=\{v_0,v_1,\dots,v_m\}=p^{\,0}(\hat{s})$.

Further, assume that
$p^{\,0}(\hat{s}')=p^{\,0}(\hat{s})=\{v_0,v_1,\dots,v_m\}$. Then
according to (\ref{def p and G action}),
$\hat{s}'=\{(v_0,g'_0),(v_1,g'_1),\dots,(v_m,g'_m)\}$, where
$g'_0,g'_1,\dots,g'_m$ are some elements of group $G$, and
$g_0+g_i=J([v_0v_i])=g'_0+g'_i$ for $i=1,\dots,m$. Set
$g=g'_0+g_0$. Then according to the above equalities $g'_i=g+g_i$
for all $i=0,1,\dots,m$ and hence $\hat{s}'=g\cdot\hat{s}$.

This proves that $p^{\,0}$ and $\lambda^0$ have the following property:
\begin{itemize}
\item[\rm(C3)] for arbitrary abstract simplices $\hat{s},\hat{s}'\in\hat{K}$
the equality $p^{\,0}(\hat{s})=p^{\,0}(\hat{s}\,')$
is equivalent to the existence of an element
$g\in{G}$ such that $g\cdot\hat{s}=\hat{s}\,'$.
\end{itemize}

It is known that $p^{\,0}$ and $\lambda^0$ may have the unique
continuation to the simplicial mapping $p:\hat{P}\to{P}$ and the
simplicial action $\lambda:G\times\hat{P}\to\hat{P}$ of group $G$
on $\hat{P}$. It also follows from (C1) -- (C3) that $p$ is a
regular covering, and $G$ is a corresponding group of covering
transformations (see, for example, \cite{Spn}).
\end{proof}

\begin{proposition}\label{covering paths and G}
Let $x=[v_0v_1]+[v_1v_2]+\dots+[v_{s-1}v_s]$ and
$y=[u_0u_1]+[u_1u_2]+\dots+[u_{t-1}u_t]$ be edge paths of the
polyhedron $P$, that run from the vertex $v_0=u_0$ to the vertex
$v_s=u_t$,
$\hat{x}=[\hat{v}_0\hat{v}_1]+[\hat{v}_1\hat{v}_2]+\dots+[\hat{v}_{s-1}\hat{v}_s]$
and
$\hat{y}=[\hat{u}_0\hat{u}_1]+[\hat{u}_1\hat{u}_2]+\dots+[\hat{u}_{t-1}\hat{u}_t]$
paths of $\hat{P}$, that cover the paths $x$ and $y$ respectively
and have the same beginning $\hat{v}_0=\hat{u}_0$. Then
$\hat{v}_s=\hat{u}_t$ if and only if $x\sim{y}$.
\end{proposition}

\begin{proof}
Let $z=[w_0w_1]+[w_1w_2]+\dots+[w_{s-1}w_s]$ be a path in the
polyhedron $P$ and $g_0\in{G}=\Bbb{Z}_2^r$. Then the unique path
$\hat{z}$ of the polyhedron $\hat{P}$, starting in the vertex
$\hat{w}_0=(w_0,g_0)$ and covering the path $z$, is defined by the
formulas
\begin{equation}\label{vertexes for covering path}
\hat{w}_i=\left(w_i,g_0+J(z_i)\right),\quad i=1,\dots,s,
\end{equation}
where $z_i=[w_0w_1]+[w_1w_2]+\dots+[w_{i-1}w_i]$, and
\begin{equation}\label{construction of covering path}
\hat{z}=[\hat{w}_0\hat{w}_1]+[\hat{w}_1\hat{w}_2]+\dots+[\hat{w}_{s-1}\hat{w}_s].
\end{equation}

Set $g_i=g_0+J(z_i)$ for $i=1,\dots,s$ and $z_0=0$. Then
$J(z_i)=J(z_{i-1})+J([w_{i-1}w_i])$ for all $i=1,\dots,s$. At the
same time $g_i=g_{i-1}+J([w_{i-1}w_i])$ and the vertices
$\hat{w}_{i-1}$ and $\hat{w}_i$ from $\hat{V}$, defined by the
formula (\ref{vertexes for covering path}), are connected by the
edge $[\hat{w}_{i-1}\hat{w}_i]\in\hat{K}$. Then in the polyhedron
$\hat{P}$ there is defined a path (\ref{construction of covering
path}) starting at the vertex $\hat{w}_0=(w_0,b_0)$. As
$p(\hat{w}_i)=p((w_i,g_0+J(z_i)))=w_i$ for all $i=0,1,\dots,s$,
then $\hat{z}$ covers the path $z$. Since $p$ is a covering then
the path $\hat{z}$ is unique.

Assume now that $\hat{v}_0=(v_0,g_0)$, where $g_0\in{G}$. By the
above proof, the equalities $p(\hat{x})=x$, $p(\hat{y})=y$ and
$\hat{v}_0=\hat{u}_0$ imply that $\hat{v}_s=(v_s,g_0+J(x))$ and
$\hat{u}_t=(u_t,g_0+J(y))$. So $\hat{v}_s=\hat{u}_t$ if and only
if $J(x)=J(y)$. According to proposition \ref{Index function
properties}, the last equality is equivalent to the homology of
the chains $x$ and $y$.
\end{proof}

\section{Minimal Cycles Searching}

Let $E(P)=K^1(P)$ be the set of edges of the polyhedron $P$,
and $L:E(P)\to\Bbb{R}$ be a non-negative function.
Using the formulas
\begin{equation}\label{weight function}
L(0)=0\mbox{ and }
L(\{a_1,\dots,a_s\})=\sum_{i=1}^s{L}(a_i).
\end{equation}
we can extend $L$ to the function $L:C_1(R)\to\Bbb(R)$. This
function is often called weight function. And for an arbitrary
chain $x\in{C}_1(P)$ the value $L(x)$ is called its weight (see,
for example, \cite{KLR}).

Let define a weight function $\hat{L}:C_1(\hat{P})\to\Bbb{R}$ assuming that
\begin{equation}\label{weight function for covering}
\hat{L}(\hat{x})=L(p(\hat{x}))
\end{equation}
for an arbitrary chain $\hat{x}\in{C}_1(\hat{P})$.

\begin{algorithm}\label{Minimal cycle with fixed vertex and index}
{\bf Searching for the minimal cycle with fixed vertex and index.}
\smallskip

\noindent {\bf Input:}

1) list $V(P)$ of vertices for polyhedron $P$;

2) lists $U(v,P)$ of vertices incident to $v$ for all vertices
$v\in{V}(P)$;

3) index function $J:C_1(P)\to\Bbb{Z}_2^r$ relative to some basis
of group $H_{n-1}(P)$;

4) weight function $L:C_1(P)\to\Bbb{R}$;

5) vector $i\in{G}=\Bbb{Z}_2^r$;

6) vertex $u\in{V}(P)$.
\smallskip

\noindent {\bf Output:}

1-chain $z\in{C}_1(P)$.
\smallskip

\noindent{\bf Algorithm Description.}
\smallskip

{\bf Step 1. Initializing cycle $z$.} Set $z:=\emptyset$.

{\bf Step 2. Initializing sets $\hat{T}\subset V(P)\times{G}$,
$\hat{P}^*\subset V(P)\times{G}$ and a mapping
$\hat{D}:V(P)\times{G}\to\Bbb{R}$.} Let $\hat{T}:=\{(u,0)\}$,
where $0$ -- null vector of space $G=\Bbb{Z}_2^r$,
$\hat{P}^*:=\emptyset$ and $\hat{D}(u,0):=0$.

{\bf Step 3. First extension of $\hat{P}^*$ and $\hat{D}$.} For
each vertex $v\in{U}(u,P)$ set $j:=J([uv])$ and add the pair
$(v,j)$ into the list $\hat{P}^*$. At the same time set
$\hat{D}(v,j):=L([uv])$, $F(v,j):=(u,0)$.

{\bf Step 4. Choosing a next element to add to $\hat{T}$.} Find
the pair $(w,k)\in(\hat{P}^*\setminus\hat{T})$ such that
$\hat{D}(w,k)=\min\limits_{(v,j)\in(\hat{P}^*\setminus\hat{T})}
\hat{D}(v,j)$.

{\bf Step 5. Stop criterion of $\hat{T}$, $\hat{P}^*$, $\hat{D}$
construction.} If $w=u$ è $k=i$, then go to step 9.

{\bf Step 6. Extension of the set $\hat{T}$.} Add the pair $(w,k)$
into the list $\hat{T}$.

{\bf Step 7. Next extension of $\hat{P}^*$ and $\hat{D}$.} For
each vertex $v\in{U}(w,P)$ set $j:=k+J([wv])$. If the pair
$(v,j)\notin\hat{P}^*$, then set
$\hat{D}(v,j):=\hat{D}(w,k)+L([wv])$, $F(v,j):=(w,k)$ and add the
pair $(v,j)$ into $\hat{P}^*$. If
$(v,j)\in(\hat{P}^*\setminus\hat{T})$ and $\hat{D}(w,k)+L([wv]) <
\hat{D}(v,j)$, then set $\hat{D}(v,j)=\hat{D}(w,k)+L([wv])$ and
$F(v,j):=(w,k)$.

{\bf Step 8. Continuation of $\hat{T}$, $\hat{P}^*$, $\hat{D}$
construction.} Go back to step 4.

{\bf Step 9. Construction of chain $z$.}

{\bf Step 9.1.} Take a pair $(v,j) = F(w,k)$ and set $z:=z+[vw]$.

{\bf Step 9.2.} If $(v,j)\ne(u,0)$, then set $(w,k)$ equal to
$(v,j)$ and go back to step 9.1.
\smallskip

\noindent{\bf End of algorithm.}
\end{algorithm}

\begin{theorem}\label{th of Minimal cycle with fixed vertex and index}
The chain $z\in{C}_1(P)$ computed by the algorithm \ref{Minimal
cycle with fixed vertex and index} has the following properties:
\begin{itemize}
\item $z\in{Z}_1(P)$;
\item $J(z)=i$;
\item $u\in{V}(z)$, where $V(z)$ is the vertex
set of the chain $z$;
\item $L(z)\le{L}(x)$ for all
cycles $x\in{Z}_1(P)$ that satisfy conditions $J(x)=i$ and
$u\in{V}(x)$.
\end{itemize}
\end{theorem}

\begin{proof}
Let $T^*$ be the result set of Dijkstra's algorithm for a
one-dimensional skeleton $\hat{P}^1$ of the polyhedron $\hat{P}$
if we choose the pair $\hat{u}=(u,0)$ as the start point, and the
pair $(u,i)$ as the end point (see, for example, \cite{AHU}).

According to the definition of the complex $\hat{K}$, the pairs
$\hat{v}=(v,j)$ and $\hat{u}=(u,0)$ in step 3, as well as the
pairs $\hat{v}=(v,j)$ and $\hat{w}=(w,k)$ in step 7 are connected
by the edges $[\hat{v}\hat{u}]\in E(\hat{P})$ and
$[\hat{v}\hat{w}]\in E(\hat{P})$ respectively. Also, according to
(\ref{weight function for covering}), we have the equalities
$\hat{L}([\hat{v}\hat{u}])=L([vu])$ in step 3 and
$\hat{L}([\hat{v}\hat{w}])=L([vw])$ in step 7. This implies that
the set $\hat{T}$ constructed by step 9 is the same that $T^*$.

Let note that step 9 is not limited to compute
$z=[v_0v_1]+\dots+[v_{q-1}v_q]$ starting at $v_0=u$ and ending at
$v_q=u$, but it also gives us the possibility to construct the
vector sequence $j_0,j_1,\dots,j_q\in\Bbb{Z}_2^r$, that will
satisfy the equalities $j_0=0$, $j_q=i$ and
$j_s=j_{s-1}+J([v_{s-1}v_s])$.

Set $\hat{v}_s=(v_s,j_s)$ for all $s=0,1,\dots,q$. Then
$[\hat{v}_{s-1}\hat{v}_s]\in{E}(\hat{P})$ for the same $s$ and
$\hat{z}=
[\hat{v}_0\hat{v}_1]+[\hat{v}_1\hat{v}_2]+\dots+[\hat{v}_{q-1}\hat{v}_q]$
is a path in the skeleton $\hat{P}^1$, starting at $\hat{u}=(u,0)$
and ending at $(u,i)$. Since it can be computed by Dijkstra's
algorithm, $\hat{L}(\hat{z})$ is not over than weight of any other
path in $\hat{P}^1$, running from $\hat{u}=(u,0)$ to $(u,i)$.

By the construction of the path $\hat{z}$ and according to
(\ref{weight function for covering}), $L(z)=\hat{L}(\hat{z})$.
Now, in the polyhedron $P$, let consider another cycle $z'$
containing the vertex $u$ and having the index $J(z')=i$.
According to proposition \ref{Index function properties},
$[z']=[x]$ in $H_1(P)$. Since $p:\hat{P}\to{P}$ is a covering,
there exists the unique path $\hat{z}'$ in $\hat{P}$, that covers
$z'$ and starts at the vertex $\hat{u}=(u,0)$. At the same time,
by statements \ref{Index function properties} and \ref{covering
paths and G} the end points of these paths $\hat{z}$ and
$\hat{z}'$ coincide. But then according to the above proof,
$L(z)=\hat{L}(\hat{z})\le\hat{L}(\hat{z}')=L(z')$.
\end{proof}

\begin{algorithm}\label{Minimal cycle from fixed gomology class}
{\bf Searching for the minimal cycle from fixed gomology class.}
\smallskip

\noindent {\bf Input:}

1) list $V(P)$ of vertices for polyhedron $P$;

2) lists $U(v,P)$ of vertices incident to $v$ for all vertices
$v\in{V}(P)$;

3) simple basis cycles $z_1^{n-1},z_2^{n-1},\dots,z_r^{n-1}$ of
homology group $H_{n-1}(P)$;

4) lists $V(z_1^{n-1}),\dots,V(z_r^{n-1})$ of vertices from cycles
$z_1^{n-1},\dots,z_r^{n-1}$ respectively;

5) index function $J:C_1(P)\to\Bbb{Z}_2^r$ relative to basis
$[z_1^{n-1}],\dots,[z_r^{n-1}]$ of group $H_{n-1}(P)$;

6) weight function $L:C_1(P)\to\Bbb{R}$;

7) cycle $x\in{Z}_1(P)$.
\smallskip

\noindent {\bf Output:}

1-cycle $z\in{Z}_1(P)$.
\smallskip

\noindent{\bf Algorithm Description.}
\smallskip

{\bf Step 1.} Set $Z:=\emptyset$.

{\bf Step 2.} Determine the vector $i=J(x)$.

{\bf Step 3.} If $i=0$, then set $z=0$ and go to step 7.

{\bf Step 4.} Find a number $k\in\{1,\dots,r\}$ such that
coordinate $i^k$ of the vector $i$ is equal 1.

{\bf Step 5.} For each vertex $v\in{V}(z_k^{n-1})$ execute steps
5.1 -- 5.3.

{\bf Step 5.1.} Using algorithm \ref{Minimal cycle with fixed
vertex and index} we find containing $v$ cycle $z_v\in Z_1(P)$
with index $J(z_v)=i$ having minimal weight $L(z_v)$ in set of all
cycles with the same properties.

{\bf Step 5.2.} Add the cycle $z_v$ into the list $Z$.

{\bf Step 5.3.} Take the next vertex $v\in{V}(z_k^{n-1})$.

{\bf Step 6.} Choose the cycle $z\in Z$ such that
$L(z)=\min\limits_{z'\in{Z}}L(z')$.

{\bf Step 7.} Quit.
\smallskip

\noindent{\bf End of algorithm.}
\end{algorithm}

\begin{theorem}\label{th of Minimal cycle from fixed gomology class}
Let $z$ be the the cycle found by the algorithm \ref{Minimal cycle
from fixed gomology class}. Then
\begin{itemize}
\item $z\sim{x}$;
\item $L(z)=\min\limits_{y\in[x]}L(y)$.
\end{itemize}
\end{theorem}

\begin{proof}
First, if $i=0$, then according to proposition \ref{Index function
properties}, cycle $x$ is homologous to zero. At the same time we
assume in step 3 that $z=0$. According to (\ref{weight function}),
$L(0)=0$. Thus, in this case $z\sim x$ and
$L(z)=\min\limits_{y\in[x]}L(y)$.

Further, let $i\ne0$. Then according to step 4
$i^k=1$ for $k\in\{1,\dots,r\}$.

Let now consider an arbitrary element $z_v$ in the list $Z$. It is
chosen in step 5.1, and according to this step $J(z_v)=i=J(x)$.
According to proposition \ref{Index function properties} it
follows that $z_v\sim x$. Since $z=z_v$ for some
$v\in{V}(z_s^{n-1})$ then $z\sim{x}$ too.

Let assume that some one-dimensional cycle $y$ of the polyhedron
$P$ belongs to the class $[x]$. Then $J(y)=J(x)=i$. Hence,
$\Ind([y],[z_k^{n-1}])=J^k(y)=1$, and therefore the cycles $y$ and
$z_k^{n-1}$ have at least one common vertex $u\in{V}(z_k^{n-1})$.
In this case, according to the selection of cycle $z_u$ in step
5.1 of algorithm \ref{Minimal cycle from fixed gomology class},
$L(z_u)\le{L}(y))$. This implies according to step 6, that
$L(z)\le{L}(z_u)\le{L}(y)$.
\end{proof}

\bibliographystyle{amsplain}

\end{document}